\documentclass{llncs}
\usepackage{amsfonts}
\usepackage{makeidx}  
\makeindex
\def\mod{\hbox{\rm{mod}}}

\begin{document}

\title{The formulas  of coefficients of sum and
product of $p $-adic integers \\ with applications to Witt vectors}

\author{Kejian Xu\inst{1}   ,\ \  Zhaopeng Dai \inst{2}  \and
Zongduo Dai \inst{3}  }

\institute{College of Mathematics,Qingdao University, China,\\
\email{kejianxu@amss.ac.cn}
\and
Institute of System Science,
Academy of Mathematics and System Science,\\
Chinese Academy of Sciences,China, \\
\email{ dzpeng@amss.ac.cn }
\and
 State Key Laboratory of Information Security,
Graduate School of Chinese Academy of Sciences,China£¬\\
\email{zongduodai@is.ac.cn}}

\maketitle              
\bigskip

\begin{figure}[b]
\rule[-2.5truemm]{5cm}{0.1truemm}\\[2mm]
{\small
\begin{tabular}{ll}
 This research is supported by National Natural Science Foundation of
China (10871106, 9060\\4011, 60473025).
\end{tabular}}
\end{figure}

\bigskip

{\bf Abstract} The explicit formulas  of operations, in particular
addition and multiplication, of $p $-adic integers are presented. As
applications of the results,  at first the explicit formulas of
operations of Witt vectors with coefficients in $\mathbb{F}_{2}$ are
given; then, through solving a problem of Browkin about the
transformation between the coefficients of a $p$-adic integer
expressed in the ordinary least residue system and the numerically
least residue system, similar formulas for Witt vectors with
coefficients in $\mathbb{F}_{3}$ are obtained.

\bigskip

\section{Introduction}

 For any two $p$-adic integers $a,b \in  \mathbb{Z} _{p}$, assume that we have the $p$-adic expansions:
$$a=a_0+a_{1}p+a_{2}p^{2}+\cdots +a_{n}p^{n}+\ldots$$
$$b=b_0+b_{1}p+b_{2}p^{2}+\cdots +b_{n}p^{n}+\ldots \ $$
$$\ \ \ \ a+b=c_0+c_{1}p+c_{2}p^{2}+\cdots +c_{n}p^{n}+\ldots$$
$$\ \ -a=d_0+d_{1}p+d_{2}p^{2}+\cdots +d_{n}p^{n}+\ldots$$
$$\ ab=e_0+e_{1}p+e_{2}p^{2}+\cdots +e_{n}p^{n}+\ldots$$
then we have the following problem.

{\bf{Problem}}  {\it For any $t$,  express $c_{t}, d_{t},e_{t}$ by
some polynomials over $\mathbb{F}_{p}$ of $a_{0},a_{1},\cdots  , \\  a_{t};b_{0}, b_{1}, \cdots,b_{t}$.}

In this paper, this problem is investigated. In section 2 and
section 3 of this paper,  we write out the polynomials for $c_{t}$
and $d_{t}$ explicitly. In section 4, we deal with the case of $ab,$
which is rather complicated, and we give an expression of $e_{t},$
which reduces the problem to the one about some kinds of partitions
of the integer $p^{t}.$

As an application, we apply the results to the operations on Witt
vectors([1]).  Let $R$ be an associative ring. The so-called Witt
vectors are vectors $(a_{0},a_{1},\cdots), a_{i}\in R$, with the
addition and the multiplication defined as follows.
$$(a_{0},a_{1},\ldots)\dot{+}(b_{0},b_{1},\ldots)=(S_{0}(a_{0},b_{0}),S_{1}(a_{0},a_{1};b_{0},b_{1}),\ldots)$$
$$\ \ \ (a_{0},a_{1},\ldots)\dot{\times}(b_{0},b_{1},\ldots)=(M_{0}(a_{0},b_{0}),M_{1}(a_{0},a_{1};b_{0},b_{1}),\ldots),$$
where $S_{n}, M_{n}$ are rather complicated polynomials in
$\mathbb{Z}[x_{0},x_{1},\ldots,x_{n};y_{0},y_{1},\ldots,y_{n}]$ and
can be uniquely but only recurrently determined by  Witt polynomials
(see [1]). Up to now it seems too involved to find patterns for
simplified forms of $S_{n}$ and $M_{n}$ for all $n$, and therefore
no explicit expressions for $S_{n}$ and $M_{n}$ are given yet. It is
well known that all Witt vectors with respect to the addition
$\dot{+}$ and the multiplication $\dot{\times}$ defined above form a
ring, called the ring of Witt vectors with coefficients in $R$
 and denoted by $\mathbf{W}(R)$.
 A similar problem is  whether the addition and the multiplication of Witt vectors can be expressed explicitly.
   From [1] it is well known that we have the canonical isomorphism
   $$\mathbf{W}(\mathbb{F}_{p})\cong \mathbb{Z}_{p},$$
   which is given by
$$(a_{0},a_{1},\ldots, a_{i},\ldots)\longmapsto \sum_{i=0}^{\infty}\tau (a_{i})p^{i},$$
where $\tau$ is the Teichm\"{u}ller  lifting. By this isomorphism,
the operations on $\mathbb{Z}_{p}$ can be transmitted to those on
$\mathbf{W}(\mathbb{F}_{p})$. But, here the elements of
$\mathbb{Z}_{p}$ are expressed with respect to the multiplicative
residue system $\tau(\mathbb{F}_{p})$, not the ordinary least
residue system $\{0,1, \ldots, p-1\}$. So,  for $p> 5$ the
operations on $\mathbb{Z}_{p}$ and hence on
$\mathbf{W}(\mathbb{F}_{p})$ do not coincide with the ordinary
operations of  $p$-adic integers.  While in the case of $p=2,$ we
have $\tau(\mathbb{F}_{2})=\{0,1\},$ that is, the two residue
systems coincide.
  Hence,
our results in the case of $p=2$ imply  that the operations on Witt
vectors in $\mathbf{W}(\mathbf{\mathbb{F}}_{2})$ can be written
explicitly. As for the case of $p=3,$ we have
$\tau(\mathbb{F}_{3})=\{-1,0,1\},$ but it is difficult to apply our
results directly to $\mathbf{W}(\mathbb{F}_{3}).$ However, in a
recent private communication, Browkin once considered the
transformation between the coefficients of a $p$-adic integer
expressed in the ordinary least residue system and the numerically
least residue system, and proposed the following problem, which
provides us a way to apply our results to
$\mathbf{W}(\mathbb{F}_{3}).$

{\bf{Browkin's problem}} \ {\it Let $p$ be an odd prime. Every
$p$-adic integer $c$ can be written in two forms:

$$c=\sum_{i}^{\infty}a_{i}p^{i}=\sum_{j}^{\infty}b_{j}p^{j},$$
where $a_{i}$ and $b_{j}$ belong respectively to the sets:

$\ \ \ \ \ \ \ \ \ \ \ \ \ \ \ \ \ \ \ \{0,1,\ldots,p-1\} \  $   and  $  \ \{0, \pm1,\pm2,\ldots,\pm\frac{p-1}{2}\}$

Obviously every $b_{j}$ is a
polynomial of $a_{0},a_{1},\ldots,a_{j}$ (and conversely). Can we
write these polynomials explicitly ?}

In section 5 of this paper, we solve Browkin's problem, that is, we
present  the required polynomials. And so, as an application, in
section 6 we can write the operations of
$\mathbf{W}(\mathbf{\mathbb{F}}_{3})$ explicitly.

\bigskip

\section{Addition}

By convention, for the empty set $\phi$, we let $\prod_{i\in
\phi}=1.$

{\bf Theorem 2.1.}  {\it Assume that
$$ A=\sum_{i=0}^{r} a_{i} p^{i},  B
=\sum_{i=0}^{r} b_{i} p^{i}, A+B =\sum_{i=0}^{r+1} c_{i} p^{i},$$
 where $a_{i} , b_{i}, c_{i} \in \{ 0,1,\ldots, p-1 \}  $ and $r\geq 1.$ Then
  $c_{0} = a_{0}+b_{0}\,( \mod \,p),$
   and  for $1\leq t\leq r+1,$
$$c_{t}=a_{t}+b_{t}+\sum_{i=0}^{t-1}\left(\sum_{k=1}^{p-1}\left(\begin{array}{c} a_{i} \\ k\\
 \end{array}
  \right)\left(\begin{array}{c} b_{i} \\ p-k\\
 \end{array}
  \right) \right)  \prod_{j=i+1} ^{t-1}\left(\begin{array}{c} a_{j}+b_{j} \\ p-1\\
 \end{array}
  \right)( \mod  p).$$}

{\bf Proof} In order to prove our result, we need the following two
lemmas.

{\bf  Lemma 2.2.} (Lucas)   {\it  If $ A=\sum_{i=0}^{r} a_{i} p^{i},
$ $ B =\sum_{i=0}^{r} b_{i} p^{i},$ $ 0 \leq a_{i} < p \, ,\,$ $ 0
\leq b_{i} < p \, ,$ then
$$ \left(\begin{array}{c} A \\ B \\
 \end{array}
  \right)    = \prod_{i=0}^{r}\left(
                                \begin{array}{c}
                                  a_{i} \\
                                   b_{i}\\
                                \end{array}
                              \right)   ( \mod \,\, p ).
$$
 In particular
$$a_{t}= \left( \begin{array}{c}
                                  A \\
                                    p^{t}\\
                                \end{array}
                              \right   )  ( \mod \,\, p ),  \ \ \  \forall \
                              t.
                              $$}

For the convenience of  readers, we include a short proof. In
$\mathbb{F}_{p}[z]$ we have
$$\sum_{t=0}^{A}\left( \begin{array}{c}
                                  A \\
                                    t\\
                                \end{array}
                              \right   ) z^{t}=(1+z)^{A}=\prod_{i=0}^{r}(1+z)^{a_{i}p^{i}}$$
                              $$=\prod_{i=0}^{r}(1+z^{p^{i}})^{a_{i}}
                              =\prod_{i=0}^{r}\sum_{j=0}^{p-1}\left( \begin{array}{c}
                                  a_{i} \\
                                    j\\
                                \end{array}
                              \right   )z^{jp^{i}} $$
$$\ \ \ \ \ \ \ \ \ \ \ \ \ =\sum_{\begin{array}{c} (j_{0},\ldots, j_{r})\\ 0\leq j_{i}\leq p-1\end{array}}\left( \begin{array}{c}
                                  a_{0} \\
                                    j_{0}\\
                                \end{array}\right)\left( \begin{array}{c}
                                  a_{1} \\
                                    j_{1}\\
                                \end{array}\right)\cdots \left( \begin{array}{c}
                                  a_{r} \\
                                    j_{r}\\
                                \end{array}\right)z^{\sum_{i=0}^{r}j_{i}p^{i}}.$$
Comparing coefficients of $z^{B}$ in both sides we get the lemma.

{\bf Lemma 2.3.} $\left( \begin{array}{c}
                                  A+B \\
                                    t\\
                                \end{array}
                              \right   ) =\sum_{\lambda+\mu=t}\left( \begin{array}{c}
                                  A\\
                                    \lambda\\
                                \end{array}
                              \right   )\left( \begin{array}{c}
                                  B \\
                                    \mu\\
                                \end{array}
                              \right   ).$

In fact, we have
$$\sum_{t}\left( \begin{array}{c}
                                  A+B \\
                                    t\\
                                \end{array}
                              \right   ) z^{t}=(1+z)^{A+B}=(1+z)^{A}(1+z)^{B}$$
                              $$\ \ \ \ \ \ \ \ \ \ \ \ \ \ =\sum_{\lambda}\left( \begin{array}{c}
                                  A\\
                                    \lambda\\
                                \end{array}
                              \right   )z^{\lambda}\sum_{\mu}\left( \begin{array}{c}
                                  B \\
                                    \mu\\
                                \end{array}
                              \right   )z^{\mu}=\sum_{t}\left (\sum_{\lambda+\mu=t}\left( \begin{array}{c}
                                  A\\
                                    \lambda\\
                                \end{array}
                              \right   )\left( \begin{array}{c}
                                  B \\
                                    \mu\\
                                \end{array}
                              \right   )\right ) z^{t}.$$
Then, the lemma follows from comparing coefficients of $z^{t}$ in
both sides.

 Now, we turn to the proof of the theorem. By the two lemmas, we have
$$c_{t}=a_{t}+b_{t}+\sum_{\lambda+\mu=p^{t}, p^{t-1}\parallel \lambda}\left(\begin{array}{c} A \\ \lambda \\
 \end{array}
  \right)\left(\begin{array}{c} B \\ \mu\\
 \end{array}
  \right)+\sum_{i=0}^{t-2}\sum_{\lambda+\mu=p^{t}, p^{i}\parallel \lambda}\left(\begin{array}{c} A \\ \lambda \\
 \end{array}
  \right)\left(\begin{array}{c} B \\ \mu\\
 \end{array}
  \right)( \mod \ p).$$
Let
$$\lambda=\lambda_{i}p^{i}+\lambda_{i+1}p^{i+1}+\ldots+\lambda_{t-1}p^{t-1},$$
where $1\leq \lambda_{i} \leq p-1, 0\leq \lambda_{j} \leq p-1$ for
$i+1\leq j \leq t-1.$ Then
$$\mu=p^{t}-\lambda=(p-\lambda_{i})p^{i}+(p-1-\lambda_{i+1})p^{i+1}+\ldots+(p-1-\lambda_{t-1})p^{t-1}. $$
Consequently, by Lucas lemma, we have in $\mathbb{F}_{p}$
$$\left(\begin{array}{c} A \\ \lambda\\
 \end{array}
  \right)=\left(\begin{array}{c} a_{i} \\ \lambda_{i}\\
 \end{array}
  \right)\prod_{j=i+1} ^{t-1}\left(\begin{array}{c} a_{j} \\ \lambda_{j}\\
 \end{array}
  \right), \ \ \left(\begin{array}{c} B \\ \mu\
 \end{array}
  \right)=\left(\begin{array}{c} b_{i} \\ p-\lambda_{i}\\
 \end{array}
  \right)\prod_{j=i+1} ^{t-1}\left(\begin{array}{c} b_{j} \\ p-1-\lambda_{j}\\
 \end{array}
  \right),$$
$$\sum_{\lambda+\mu=p^{t}, p^{t-1}\parallel \lambda}\left(\begin{array}{c} A \\ \lambda \\
 \end{array}
  \right)\left(\begin{array}{c} B \\ \mu\\
 \end{array}
  \right)=\sum_{i=1}^{p-1}\left(\begin{array}{c} a_{t-1} \\ i \\
 \end{array}
  \right)\left(\begin{array}{c} b_{t-1} \\ p-i \\
 \end{array}
  \right).$$
Therefore
$$\sum_{\lambda+\mu=p^{t}, p^{i}\parallel \lambda}\left(\begin{array}{c} A \\ \lambda \\
 \end{array}
  \right)\left(\begin{array}{c} B \\ \mu\\
 \end{array}
  \right)=\sum_{\lambda_{i}=1}^{p-1}\sum_{\lambda_{i+1}=0}^{p-1}\cdots
  \sum_{\lambda_{t-1}=0}^{p-1}\left(\begin{array}{c} a_{i} \\ \lambda_{i}\\
 \end{array}
  \right)\left(\begin{array}{c} b_{i} \\ p-\lambda_{i}\\
 \end{array}
  \right)\prod_{j=i+1} ^{t-1}\left(\begin{array}{c} a_{j} \\ \lambda_{j}\\
 \end{array}
  \right)\left(\begin{array}{c} b_{j} \\ p-1-\lambda_{j}\\
 \end{array}
  \right).$$
$$=\sum_{\lambda_{i}=1}^{p-1}\left(\begin{array}{c} a_{i} \\ \lambda_{i}\\
 \end{array}
  \right)\left(\begin{array}{c} b_{i} \\ p-\lambda_{i}\\
 \end{array}
  \right) \sum_{\lambda_{i+1}=0}^{p-1}\left(\begin{array}{c} a_{i+1} \\ \lambda_{i+1}\\
 \end{array}
  \right)\left(\begin{array}{c} b_{i+1} \\ p-1-\lambda_{i+1}\\
 \end{array}
  \right) \cdots \sum_{\lambda_{t-1}=0}^{p-1}\left(\begin{array}{c} a_{t-1} \\ \lambda_{t-1}\\
 \end{array}
  \right)\left(\begin{array}{c} b_{t-1} \\ p-1-\lambda_{t-1}\\
 \end{array}
  \right)$$
To all of these sums but the first we apply Lemma 2.3
  and we get
  $$\sum_{\lambda_{i}=1}^{p-1}\left(\begin{array}{c} a_{i} \\ \lambda_{i}\\
 \end{array}
  \right)\left(\begin{array}{c} b_{i} \\ p-\lambda_{i}\\
 \end{array}
  \right)\cdot \prod_{j=i+1}^{t-1}\left(\begin{array}{c} a_{j}+ b_{j} \\ p-1\\
 \end{array}
  \right).$$
Therefore
$$c_{t}=a_{t}+b_{t}+\sum_{k=1}^{p-1}\left(\begin{array}{c} a_{t-1} \\ k \\
 \end{array}
  \right)\left(\begin{array}{c} b_{t-1} \\ p-k \\
 \end{array}
  \right)+\sum_{i=0}^{t-2}\left(\sum_{k=1}^{p-1}\left(\begin{array}{c} a_{i} \\ k\\
 \end{array}
  \right)\left(\begin{array}{c} b_{i} \\ p-k\\
 \end{array}
  \right)\right)\prod_{j=i+1} ^{t-1}\left(\begin{array}{c} a_{j}+b_{j} \\ p-1\\
 \end{array}
  \right)$$
$$=a_{t}+b_{t}+\sum_{i=0}^{t-1}\left(\sum_{k=1}^{p-1}\left(\begin{array}{c} a_{i} \\ k\\
 \end{array}
  \right)\left(\begin{array}{c} b_{i} \\ p-k\\
 \end{array}
  \right) \right)  \prod_{j=i+1} ^{t-1}\left(\begin{array}{c} a_{j}+b_{j} \\ p-1\\
 \end{array}
  \right)( \mod \ p).\ \ \ \ \ \ \ \ \ \ \ \ $$ $ \Box $

{\bf Corollary  2.4. } {\it Assume that
$$a=\sum_{i=0}^{\infty}a_{i}p^{i}, b=\sum_{i=0}^{\infty}b_{i}p^{i}, a+b=\sum_{i=0}^{\infty}c_{i}p^{i}  \in
\mathbb{Z}_{p},
$$ with $a_{i}, b_{i}, c_{i}\in \{0,1,\ldots,p-1\}.$ Then
 $c_{0} = a_{0}+b_{0} (\mod\,p)$,
  and  for $t\geq 1,$
$$c_{t}=a_{t}+b_{t}+\sum_{i=0}^{t-1}\left(\sum_{j=1}^{p-1}\left(\begin{array}{c} a_{i} \\ j\\
 \end{array}
  \right)\left(\begin{array}{c} b_{i} \\ p-j\\
 \end{array}
  \right)\right )\prod_{j=i+1} ^{t-1}\left(\begin{array}{c} a_{j}+b_{j} \\ p-1\\
 \end{array}
  \right)( \mod \ p).$$
 In particular, if $p=2,$ then we have $c_{0} = a_{0}+b_{0}( \mod \,2)$,  and for $t\geq 1,$
$$c_{t} = a_{t}+b_{t}+\sum_{i=0}^{t-1}a_{i}b_{i}\prod_{j=i+1}^{t-1}(a_{j}+b_{j})( \mod \ 2).$$
 }
$ \Box $

{\bf Corollary 2.5.}  {\it Assume that
$a=\sum_{i=0}^{\infty}a_{i}2^{i} \in \mathbb{Z}_{2}$ and $n\geq1.$}

(i) {\it If $2^{n}a=\sum_{i=0}^{\infty}c_{i}2^{i}  \in
\mathbb{Z}_{2},$ then $c_{t}=0, 0\leq t < n$ and
$c_{t}=a_{t-n}(\mod\,2)$  for $t\geq n.$ }

(ii) {\it If $(2^{n}+1)a=\sum_{i=0}^{\infty}c_{i}2^{i} \in
\mathbb{Z}_{2},$ then $ c_{t}=a_{t}, 0\leq t \leq n-1,
c_{n}=a_{n}+a_{0} ( \mod \ 2)$ and for $t\geq n+1,$
$$c_{t}=a_{t}+a_{t-n}+\sum_{i=n}^{t-1}a_{i}a_{i-n}\prod_{j=i+1}^{t-1}(a_{j}+a_{j-n})( \mod \
2).$$} $\Box $

{\bf Corollary 2.6.}  {\it Assume that
$a=\sum_{i=0}^{\infty}a_{i}3^{i} \in \mathbb{Z}_{3}$ and $n\geq1.$
 If $2a=\sum_{i=0}^{\infty}c_{i}3^{i} \in
\mathbb{Z}_{3},$ then $c_{0}=-a_{0} ( \mod \,3) $ and for $t\geq
1,$
$$c_{t}=-a_{t}+\sum_{i=0}^{t-1}
  a_{i}(1-a_{i})  \prod_{j=i+1} ^{t-1}a_{j}(2a_{j}-1)( \mod \
  3).$$}
$\Box $

\bigskip

\section{Minus}

{\bf Theorem 3.1.}  \ {\it Let $A=\sum_{i=0}^{r} a_{i} p^{i}.$
Assume that
$$ -A=\sum_{i=0}^{r} d_{i} p^{i}\, ( \mod  \, p^{r+1}),$$
  where
$d_{i}\in \{0,1,\ldots,p-1 \}. $ Then $d_{0}=-a_{0} (\mod\,p)$ and
for  $1\leq t \leq r$
$$d_{t}=-a_{t}-1+ \prod_{i=0}^{t-1}(1-a_{i}^{p-1})( \mod  \,
p).$$}

{\bf Proof  }  \ \ Clearly, we can assume that $A\neq 0.$ In this
case, there exists an $s$ such that $a_{s}\neq 0 $ but $a_{i}=0 $
for $ i < s .$ This implies  that
$$d_{t}=\{\begin{array}{c} -a_{t} \ ( \mod \, p) ,\ \ \ \   \mbox{if}  \ t\leq s; \\ -a_{t}-1  (\mod \, p), \mbox{if} \  t> s,  \end{array}$$
which is equivalent to
$$d_{t}=\{\begin{array}{c} -a_{t} \ (\mod \, p), \ \ \ \ \   \mbox{if} \ (a_{0}, a_{1},\ldots, a_{t-1})=(0, 0, \ldots, 0); \\ -a_{t}-1 \ (\mod \, p),  \mbox{if} \  (a_{0}, a_{1},\ldots, a_{t-1})\neq(0, 0, \ldots, 0).  \end{array}$$
Take $f(a_{0},a_{1},\ldots,
a_{t-1})=-1+\prod_{i=0}^{t-1}(1-a_{i}^{p-1})(\mod  \, p).$
Clearly
$$f(a_{0},a_{1},\ldots, a_{t-1})=\{ \begin{array}{c} 0 \ (\mod \, p), \ \  \mbox{if} \ (a_{0}, a_{1},\ldots, a_{t-1})=(0, 0, \ldots, 0); \\
-1 (\mod \, p),   \mbox{if} \ (a_{0}, a_{1},\ldots,
a_{t-1})\neq(0, 0, \ldots, 0)\end{array}.$$ Therefore,
$$d_{t}=-a_{t}+f(a_{0},a_{1},\ldots, a_{t-1})=-a_{t}-1+ \prod_{i=0}^{t-1}(1-a_{i}^{p-1})(\mod  \, p).$$ $\Box$

{\bf Corollary 3.2.} \  {\it Assume that
$$a=\sum_{i=0}^{\infty}a_{i}p^{i}, -a=\sum_{i=0}^{\infty}d_{i}p^{i} \in
\mathbb{Z}_{p},
$$ with $a_{i}, d_{i}\in \{0,1,\ldots,p-1\}.$ Then $d_{0}=-a_{0}
(\mod \,p)$  and for $t\geq1$
$$d_{t}=-a_{t}-1+ \prod_{i=0}^{t-1}(1-a_{i}^{p-1})(\mod  \,
p).$$ If $p=2,$ then $d_{0}=a_{0},$ and for $t\geq1,$
$$ d_{t} = a_{t} +  1+ \prod_{i=0}^{t-1}(1+a_{i})
(\mod \,2).$$ } $\Box$

 {\bf Remark 3.4} The problems
considered in this section and in Corollary 2.5 and 2.6 were
suggested to us by Browkin.

\bigskip

\section{Multiplication}

{\bf 4.1. Fundamental lemma}

{\bf 4.1.1. Fundamental polynomials}  Let
$$\mathbb{K}=\{\underline{k}=(k_{1},\ldots,k_{l},\ldots,k_{p-1}):  k_{l}\geq 0, 0 \leq \sum _{l=1}^{p-1}k_{l}\leq p-1\}.$$
Clearly $\underline{0}=(0,\ldots,0)\in \mathbb{K}.$ Let
$$\mathbb{K}^{(r+1)^{2}}=\underbrace{\mathbb{K}\times
\mathbb{K}\times\cdots \times\mathbb{K}}_{(r+1)^{2}},$$ and write
$\underline{\underline{0}}=(\underline{0},\ldots,\underline{0})\in
\mathbb{K}^{(r+1)^{2}} .$

 For any
$\underline{k}=(k_{1},\ldots,k_{l},\ldots,k_{p-1})\in \mathbb{K},
\underline{k}\neq \underline{0},$ define
$$\pi_{\underline{k}} (x,y) =\frac{ y(y-1)\cdots (y-\sum_{l=1}^{p-1} k_{l}  +1 )}{ k_{1} !\cdots k_{p-1}!}\prod_{l=1}^{p-1}\left( \frac{ x(x-1)\cdots (x-l  +1 )}{l!}\right ) ^{k_{l}} \
 ( \mod  \, p),$$
and for $\underline{k}=\underline{0},$  define $\pi_{\underline{k}}
(x,y)=1$.

Let $\mathbf{I}=\{(i,j) : 0\leq i,j \leq r\},$ and let
$\underline{x}=(x_{0},\ldots, x_{r}), \underline{y}=(y_{0},\ldots,
y_{r}).$ Then for
$\underline{\underline{k}}=(\ldots,\underline{k}_{i,j},\ldots)\in
\mathbb{K}^{(r+1)^{2}}$ with $\underline{k}_{i,j}=(k_{i,j,1},\ldots,
k_{i,j,p-1}),$ we define the function
$$\pi_{\underline{\underline{k}}}(\underline{x},\underline{y})=\prod_{(i,j)\in \mathbf{I}}\pi_{\underline{k}_{i,j}} (x_{i},y_{j}),$$
 and the norm
$$\|\underline{\underline{k}} \|=\sum_{(i,j)\in \mathbf{I}}\left(\sum _{l=1}^{p-1}lk_{i,j,l}\right )p^{i+j}.$$
Clearly,
$\pi_{\underline{\underline{k}}}(\underline{x},\underline{y})$ is a
 polynomial in $x_{0},\ldots, x_{r}; y_{0},\ldots, y_{r}.$

{\bf Lemma 4.1.} {\it Assume that $\underline{0}\neq
\underline{k}\in \mathbb{K}.$ Let $0\leq a \leq p-1, 0\leq b \leq
p-1.$ Then we have $\pi_{\underline{k}}(a,b)=0,$ if one of the
following cases occurs.

(i) $ab=0;$

(ii) there exists an $l,$ such that $l>a$ and $k_{l}>0;$

(iii) $\sum_{l=1}^{p-1}k_{l}>b.$}

{\bf Proof} It can be checked directly. $\Box$

{\bf Lemma 4.2.} {\it Assume that $\underline{\underline{0}}\neq
\underline{\underline{k}}=(\ldots,\underline{k}_{i,j},\ldots)\in
\mathbb{K}^{(r+1)^{2}}.$ Let
$\underline{a}=(a_{0},a_{1},\ldots,a_{r})$ and
$\underline{b}=(b_{0},b_{1},\ldots,b_{r}).$ Then we have
$$\pi_{\underline{\underline{k}}}(\underline{a},\underline{b})=0,$$
if one of the following cases occurs.

(i) there exists $(i,j)\in \mathbf{I}$ such that $a_{i}b_{j}=0$ and
$\underline{k}_{i,j}\neq \underline{0};$

(ii) there exist $(i,j)\in \mathbf{I}, l> a_{i}$, such that
$k_{i,j,l}>0;$

(iii) there exists $(i,j)\in \mathbf{I},$ such that
$\sum_{l=1}^{p-1}k_{i,j,l}> b_{j}.$}

{\bf Proof} It follows from Lemma 4.1. $\Box$

{\bf 4.2.2. Fundamental lemma}

{\bf  Lemma 4.3.} {\it Assume that
$$A=\sum_{i=0}^{r}a_{i}p^{i}, \ \ B=\sum_{i=0}^{r}b_{i}p^{i}, \ \  AB=\sum_{i=0}^{2r+1}e_{i}p^{i}.$$
Then $e_{0} = a_{0}b_{0}(\mod\,p)$  and for  $1\leq t \leq 2r+1,$
$$  e_{t} = \sum_{\begin{array}{c} \underline{\underline{k}}\in \mathbb{K}^{(t+1)^{2}}\\
                  \| \underline{\underline{k}}\| = p^{t}
                \end{array} } \pi _{\underline{\underline{k}}}( \underline{a},\underline{b})\ \ ( \mod  \,
                p),
$$
 where $\underline{a}=(a_{0},a_{1},\ldots,a_{t})$ and}
$\underline{b}=(b_{0},b_{1},\ldots,b_{t}).$

{\bf Proof}  Define
$$ \mathbf{I}(\underline{a}, \underline{b})= \{ (i,j) \in \mathbf{I} :  0\leq i,j \leq t, a_{i}b_{j}\neq0\}.$$
For any integers $0< a,b < p$, define the subset of $\mathbb{K}:$
$$\mathbb{K}(a,b)=\{\underline{k}=(k_{1},\ldots,k_{l},\ldots,k_{a},0,\ldots,0)\in \mathbb{K} :  k_{l}\geq 0, 1\leq \sum _{l=1}^{a}k_{l}\leq b\}.$$
Note that $\underline{0}\notin \mathbb{K}(a,b).$ We will denote
$\underline{k}=(k_{1},\ldots,k_{a},0,\ldots,0)$ simply by
$(k_{1},\ldots,k_{a}).$ Then, for
$\underline{k}=(k_{1},\ldots,k_{a})\in \mathbb{K}(a,b),$ clearly we
have
$$\pi_{\underline{k}} (a,b)= \left(
\begin{array}{c} b \\ \underline{k} \end{array}\right )\prod_{l=1}^{a}\left(
\begin{array}{c} a \\ l \end{array}\right ) ^{k_{l}} \
 ( \mod  \, p),$$
where
$$\left(
\begin{array}{c} b \\ \underline{k} \end{array}\right )=\frac{ b!}{k_{1}!\cdots k_{a} !(b-\sum_{l=1}^{a} k_{l} )!}.$$

For $\phi \neq S\subseteq \mathbf{I}(\underline{a}, \underline{b}),$
define the subset  of $\mathbb{K}^{(t+1)^{2}}$:
$$\mathbb{K}_{S}(\underline{a},\underline{b})=\{(\ldots,\underline{k}_{i,j},\ldots)\in \mathbb{K}^{(t+1)^{2}} :  \underline{k}_{i,j}\in \mathbb{K}(a_{i},b_{j}), (i,j)\in S; \, \underline{k_{i,j}}=\underline{0}, (i,j)\notin S \}.$$
If
$\underline{\underline{k}}=(\ldots,\underline{k}_{i,j},\ldots)\in\mathbb{K}_{S}(\underline{a},\underline{b})
$ with $\underline{k_{i,j}}=(k_{i,j,1},k_{i,j,2},\ldots, k_{i,j,
a_{i}})\in \mathbb{K}(a_{i},b_{j}),$ then it is easy to show that
$$\pi_{\underline{\underline{k}}}(\underline{a}, \underline{b})=\prod_{(i,j)\in S}\pi_{\underline{k}_{i,j}}(a_{i}, b_{j})
 ( \mod  \, p).$$
and
$$\|\underline{\underline{k}} \|=\sum_{(i,j)\in S}\left(\sum _{l=1}^{a_{i}}lk_{i,j,l}\right)p^{i+j}.$$

Now, we have
$$ \sum_{ 0 \leq \lambda \leq AB} \left(
                                   \begin{array}{c}
                                     AB \\
                                     \lambda \\
                                   \end{array}
                                 \right)z^{\lambda } = ( 1 + z)^{AB} =  \prod_{\begin{array}{c}
                                                                                    0 \leq i \leq t  \\
                                                                                     a_{i}\neq 0
                                                                                  \end{array}
                                 }(1+z^{p^{i}})^{a_{i}B}
$$
$$  = \prod_{(i,j)\in \mathbf{I}(\underline{a}, \underline{b})
                                 }\left(1+\sum_{l=1}^{a_{i}}\left( \begin{array}{c} a_{i} \\ l \end{array}\right )z^{lp^{i+j}}\right)    ^{b_{j}}  \ \ \ \ \ \ \ \  \ \ \ \ \ \  $$
$$\ \ \ \ \ \ \ \  \ \ \ \ \ \ \  =\prod_{(i,j)\in \mathbf{I}(\underline{a}, \underline{b})} \left( 1+\sum
_{\underline{k} \in \mathbb{K}(a_{i},b_{j}) } \left(
\begin{array}{c} b_{j} \\ \underline{k} \end{array}\right )\prod_{l=1}^{a_{i}}\left(
\begin{array}{c} a_{i} \\ l \end{array}\right ) ^{k_{l}}
z^{\sum_{l=1}^{a_{i}}lk_{l}p^{i+j}}   \right ) $$
$$ \ \ \ \   = \prod_{(i,j)\in \mathbf{I}(\underline{a}, \underline{b})} \left( 1+\sum _{ \underline{k} \in \mathbb{K}(a_{i},b_{j})}
 \pi_{\underline{k}}(a_{i},b_{j})    z^{\sum_{l=1}^{a_{i}}lk_{l}p^{i+j}} \right ) $$
$$= 1 + \sum_{\phi\neq S\subseteq \mathbf{I}(\underline{a}, \underline{b})}\sum_{
                \underline{\underline{k}} = ( \cdots , \underline{k}_{i,j},\cdots )
                \in \mathbb{K}_{S} (\underline{a},\underline{b}) } \prod_{(i,j) \in S} \pi _{\underline{k}_{i,j}} (a_{i}, b_{j})\cdot z^{\sum_{(i,j) \in S}(\sum_{l=1}^{a_{i}}lk_{i,j,l})p^{i+j}}$$
$$   = 1 + \sum_{\phi\neq S \subseteq \mathbf{I}(\underline{a}, \underline{b})}\sum_{\begin{array}{c}
                  \underline{\underline{k}} \in \mathbb{K}_{S}(\underline{a},\underline{b} )
\end{array} }\pi_{ \underline{\underline{k}} }(\underline{a},\underline{b} )z^{\|  \underline{\underline{k}} \|}( \mod  \, p). \ \ \ \ \ \ \ \  \ \ \ \ \ \ \ \ \ \ \ \  \ \ \ \ \ \ \ \ \ \ \ \ \ \  \ \  $$
Comparing the coefficients of both sides and letting
$\lambda=p^{t}$, then from Lucas lemma, we have
$$e_{t}=\left ( \begin{array}{c}AB \\
p^{t} \\
\end{array}
\right)=\sum_{ \phi \neq S \in \mathbf{I}(\underline{a},
\underline{b})}\sum_{\begin{array}{c}
                  \underline{\underline{k}} \in \mathbb{K}_{S}(\underline{a},\underline{b} )\\
                   \|  \underline{\underline{k}} \|=p^{t}
                 \end{array} }\pi_{\underline{\underline{k}}  }(\underline{a},\underline{b} )=\sum_{\begin{array}{c} \underline{\underline{k}}\in \mathbb{K}^{(t+1)^{2}}\\
                  \| \underline{\underline{k}}\| = p^{t}
                \end{array} } \pi _{\underline{\underline{k}}}( \underline{a},\underline{b})\ \ ( \mod  \,
                p).$$
The last step follows from Lemma 4.2. $\Box$

{\bf 4.2. Multiplication formula}

{\bf 4.2.1. $T_{p}$-partitions} \  Now we shall give a simpler
formula for $e_{t}.$ Let $\mathbb{K}^{*}=\mathbb{K}\backslash \{0\}$
and $K:=| \mathbb{K}^{*}|.$ Then $|\mathbb{K}| =K+1$ and we can
write the elements of $\mathbb{K}$ as $\underline{k}(j), 0\leq j
\leq K,$ in particular, let $\underline{k}(0)=\underline{0}$ for
convenience. So
$$\mathbb{K}^{*}=\{\underline{k}(j): 1\leq j \leq K\}.$$
For $\underline{k}=(k_{1},\ldots,k_{l},\ldots,k_{p-1})\in
\mathbb{K},$ define
$$w(\underline{k})=\sum_{j=1}^{p-1}jk_{j}.$$
In the following, we fix the vector:
$$\underline{w}:=(w(\underline{k}(1)),w(\underline{k}(2)),\ldots, w(\underline{k}(K))).$$

For $\underline{l}=(l_{1}, l_{2}, \ldots, l_{K})\in \mathbb{N}^{K}$
(the cartesian  product of  $\mathbb{N},$ the set of non-negative
integers), the size of $\underline{l}$ is defined as
$$|
\underline{l}|=\sum_{j=1}^{K}l_{j},$$
 and the inner product of
$\underline{w}$ and $\underline{l}$ is defined as
$$\underline{w}\cdot \underline{l}=\sum_{j=1}^{K}w(\underline{k}(j))l_{j}.$$
For an integer $n\geq 0,$ a $T_{p}$-partition of $n$ is defined as
$$n=\sum_{j=0}^{t}(\underline{w}\cdot {\underline{l}}_{j})p^{j}, \  {\underline{l}}_{j} \in \mathbb{N}^{K}, 0\leq | \underline{l}_{j}| \leq 1+j.$$
This partition is also written as
$$\underline{\underline{l}}=(\underline{l}_{0},\ldots,\underline{l}_{m},\ldots,\underline{l}_{t}), 0\leq | \underline{l}_{m}| \leq 1+m.$$
We will use the symbol $\mathbf{L}_{p}(t)$ to denote the set of all
possible $T_{p}$-partitions of $p^{t},$ that is,
$$\mathbf{L}_{p}(t)=\{\underline{\underline{l}}=(\underline{l}_{0},\ldots,\underline{l}_{m},\ldots,\underline{l}_{t}): \sum_{j=0}^{t}(\underline{w}\cdot {\underline{l}}_{j})p^{j}=p^{t}, 0\leq | \underline{l}_{m}| \leq 1+m\}.$$

If $p=2,$ then $K=1$ and $\underline{l}_{m}$ is only a non-negative
integer, so we can write $\underline{l}_{m}=l_{m}.$ Clearly
$l_{0}=0.$ Hence, for $p=2,$ we have
$$\mathbf{L}_{2}(t)=\{\underline{\underline{l}}=(l_{1},\ldots,l_{k},\ldots, l_{t}): \sum_{ k=1}^{t} l_{k} 2^{k} =2^{t}, 0 \leq l_{k} \leq
k+1 \}.$$

If $p=3,$ then $K=5$ and we have
$$\mathbb{K}^{*}=\{\underline{k}(1)=(1,0), \underline{k}(2)=(0,1),\underline{k}(3)=(2,0), \underline{k}(4)=(1,1), \underline{k}(5)=(0,2)\},$$
and therefore $\underline{w}=(1,2,2,3,4).$ Hence, for $p=3,$ we have
$$\mathbf{L}_{3}(t)=\{\underline{\underline{l}}=(\underline{l}_{0},\ldots,\underline{l}_{k},\ldots,\underline{l}_{t}):
 \sum_{k=0}^{t}(l_{k1}+2l_{k2}+2l_{k3}+3l_{k4}+4l_{k5})3^{k}=3^{t},$$
$$ \ \ \ \ \ \ \ \ \ \ \   0\leq | \underline{l}_{k}| \leq 1+k\},$$
where $\underline{l}_{k}=(l_{k1},l_{k2},l_{k3},l_{k4},l_{k5}), 0\leq
k \leq t.$

\newpage

{\bf 4.2.2 Partitions of $\mathbf{I(m)}$ and symmetric polynomials}

\ Let $\mathbf{I}(m)=\{i: 0\leq i \leq m\}, \ 0\leq m \leq t.$ For
$\underline{l}=(l_{1}, \ldots, l_{j},\ldots, l_{K})\in
\mathbb{N}^{K}$ with $ | \underline{l}| \leq 1+m,$ we call
$\underline{S}=(S_{1},\ldots, S_{j},\ldots, S_{K})$ an
$\underline{l}$-partition of $\mathbf{I}(m),$ if it satisfies
$$S_{j}\subseteq \mathbf{I}(m),\ |S_{j}|=l_{j}, $$
$$\ S_{j}\cap S_{j^{\prime}}=\phi, \ \forall\ j\neq j^{\prime}, 1\leq
j,j^{\prime}\leq K.$$ The set of all possible
$\underline{l}$-partitions of $\mathbf{I}(m)$ is denoted by
$\mathbf{I}(m,\underline{l}),$ that is,
$$\mathbf{I}(m,\underline{l})=\{(S_{1}, S_{2},\ldots, S_{K}): S_{j}\subseteq \mathbf{I}(m),\ |S_{j}|=l_{j}, \ S_{j}\cap S_{j^{\prime}}=\phi, \ $$
$$\ \ \ \ \ \ \ \ \ \ \ \ \ \ \ \ \ \ \ \ \ \ \forall\ j\neq j^{\prime}, 1\leq
j,j^{\prime}\leq K\}.$$
 Defining
$l_{0}:=1+m-\sum_{j=1}^{K}l_{j},$ we get
$$|\mathbf{I}(m,\underline{l})|=\frac{(1+m)!}{l_{0}!l_{1}!\cdots l_{K}!}$$

For a given integer $m, 0\leq m\leq t,$ and
$\underline{l}=(l_{1},\ldots, l_{j}, \ldots, l_{K})\in
\mathbb{N}^{K}$ with $|\underline{l}| \leq 1+m,$ define the function
$$\tau_{\underline{l}}(x_{0},\ldots,x_{m}; y_{0},\ldots,y_{m})=\sum_{\underline{S}=(S_{1},\ldots, S_{j},\ldots, S_{K})\in \mathbf{I}(m, \underline{l})}\prod_{j=1}^{K}\prod_{i\in S_{j}}\pi_{\underline{k}(j)}(x_{i},y_{m-i}).$$
Clearly, $\tau_{\underline{l}}(x_{0},\ldots,x_{m};
y_{0},\ldots,y_{m})$ is a polynomial which is symmetric with respect
to the pairs $\{(x_{i}, y_{m-i}): 0\leq i\leq m\}$, that is, it is
invariant under the permutations of the pairs.

When $p=2,$ we have $K=1, \mathbb{K}=\{0,1\}$ and hence
$\underline{k}(1)=1$ as well as $l:=l_{1}=\underline{l}.$ So we have
$$\tau_{\underline{l}}(x_{0},\ldots,x_{m}; y_{0},\ldots,y_{m})=\sum_{0\leq i_{1}< \cdots < i_{l}\leq m}\prod _{k=1}^{l}x_{i_{k}}y_{m-i_{k}}=\tau _{l
}(x_{0}y_{m},x_{1}y_{m-1},\cdots ,x_{m}y_{0}),$$ where
 $\tau _{l}(X_{0},X_{1},\cdots ,X_{m})$
denote the $l$-th elementary symmetric polynomial of $X_{0},X_{1}, \cdots ,X_{m}.$

When $p=3,$ we have the ordered set $\mathbb{K}^{*}=\{(1,0),
(0,1),(2,0), (1,1), (0,2)\}.$ It is easy to check that when $x_{i},
y_{j}\in \mathbb{F}_{3},$ as polynomial functions we have
$$\tau_{\underline{l}}(x_{0},\ldots,x_{m}; y_{0},\ldots,y_{m})=\sum_{\underline{S}=(S_{1}, \ldots , S_{5})\in
\mathbf{I}(m,
\underline{l})}f_{\underline{S}}(x_{0},x_{1},\ldots,x_{m};y_{0},y_{1},\ldots,y_{m}),$$
 where
 $$f_{\underline{S}}(x_{0},x_{1},\ldots,x_{m};y_{0},y_{1},\ldots,y_{m})=\prod_{ i_{1}\in
S_{1}}x_{i_{1}}y_{m-i_{1}}\prod_{ i_{2}\in S_{2}}
x_{i_{2}}(1-x_{i_{2}})y_{m-i_{2}}$$
$$ \ \ \ \ \ \ \ \ \ \ \ \ \ \ \ \ \ \ \ \ \ \ \ \ \ \ \ \ \ \ \ \ \ \ \ \cdot\prod_{ i_{3}\in
S_{3}}x_{i_{3}}^{2}y_{m-i_{3}}(1-y_{m-i_{3}})\prod_{ i\in S_{4}\cup
S_{5}}x_{i}(1-x_{i})y_{m-i}(y_{m-i}-1).$$

{\bf 4.2.3.  Multiplication formula}

{\bf Theorem 4.4.} {\it Assume that
$$A=\sum_{i=0}^{r}a_{i}p^{i}, \ \ B=\sum_{i=0}^{r}b_{i}p^{i}, \ \  AB=\sum_{i=0}^{2r+1}e_{i}p^{i}.$$
Then $e_{0} = a_{0}b_{0} (\mod\,p)$  and for   $1\leq t \leq 2r+1,$
$$  e_{t} = \sum_{\underline{\underline{l}}=(\underline{l}_{0},\ldots,\underline{l}_{k},\ldots,\underline{l}_{t} )\in \mathbf{L}_{p}(t)}\prod_{k=0}^{t}\tau_{\underline{l}_{k}}(a_{0},\ldots,a_{k}; b_{0},\ldots, b_{k})\ \ ( \mod   \,
                p).
$$}

{\bf Proof} For $\underline{\underline{k}}=( \cdots ,
\underline{k}_{i,j},\cdots )\in \mathbb{K}^{(t+1)^{2}},$ let
$$\underline{\underline{S}}(\underline{\underline{k}})=(\underline{S}_{0},\ldots,\underline{S}_{m},\ldots,\underline{S}_{t}), \ \underline{S}_{m}=(S_{m,1},\ldots, S_{m,j},\ldots,S_{m,K}),$$
$$\underline{\underline{l}}(\underline{\underline{k}})=(\underline{l}_{0},\ldots,\underline{l}_{m},\ldots, \underline{l}_{t}), \ \underline{l}_{m}=(l_{m,1},\ldots, l_{m,j},\ldots,l_{m,K}),$$
where
$$S_{m,j}=\{i: 0\leq i \leq m, \underline{k}_{i,m-i}=\underline{k}(j)\}, \ \ |S_{m,j}|=l_{m,j}.$$
Clearly, we have
$$S_{m,j}\subseteq \mathbf{I}(m), \ S_{m,j}\cap S_{m,j^{\prime}}=\phi, \ \forall\ j\neq j^{\prime},$$
and
$$|\underline{l}_{m}|=\sum_{j=1}^{K}l_{m,j}\leq 1+m.$$
So $\underline{S}_{m}\in \mathbf{I}(m, \underline{l}_{m}),$ and
therefore
$$\underline{\underline{S}}(\underline{\underline{k}})\in \mathbf{I}(0, \underline{l}_{0})\times\mathbf{I}(1, \underline{l}_{1})\times \cdots \times \mathbf{I}(t, \underline{l}_{t}). $$

We need the following two lemmas.

{\bf Lemma 4.5.} {\it $\|\underline{\underline{k}} \|=p^{t}$ if and
only if} $\underline{\underline{l}}(\underline{\underline{k}})\in
\mathbf{L}_{p}(t).$

In fact, noting that $w(\underline{0})=0,$ we have
$$\|\underline{\underline{k}}\|=\sum_{0\leq i,j \leq t}w(\underline{k}_{i,j})p^{i+j}=\sum_{0 \leq m\leq t}\left(\sum_{0\leq i\leq m}w(\underline{k}_{i,m-i})\right)p^{m}$$
$$=\sum_{0 \leq m\leq t}\left(\sum_{\begin{array}{c} 0\leq i\leq m, \underline{k}_{i,m-i}\neq \underline{0}\end{array}}w(\underline{k}_{i,m-i})\right)p^{m}$$
$$\ =\sum_{0 \leq m\leq t}\left(\sum_{1\leq j\leq K}\sum_{i\in S_{m,j}}w(\underline{k}(j))\right)p^{m} \ \ \ \ \ \ \ \  \ \ \ \ \ \  \ \   $$
$$\ \ \ \ \ \ \ \  \  =\sum_{0 \leq m\leq t}\left(\sum_{1\leq j\leq K}l_{m,j}w(\underline{k}(j))\right)p^{m}=\sum_{0 \leq m\leq t}(\underline{w}\cdot\underline{l}_{m})p^{m},$$
as required.

{\bf Lemma 4.6.} {\it For a fixed
$(\underline{l}_{0},\ldots,\underline{l}_{m},\ldots,\underline{l}_{t})\in
\mathbf{L}_{p}(t),$ we have the bijection:}
$$\{\underline{\underline{k}}\in \mathbb{K}^{(1+t)^{2}}: \underline{\underline{l}}(\underline{\underline{k}})=(\underline{l}_{0},\ldots,\underline{l}_{m},\ldots,\underline{l}_{t})\}\longrightarrow \mathbf{I}(0,\underline{l}_{0})\times\cdots \times \mathbf{I}(t,\underline{l}_{t})$$
$$\underline{\underline{k}}\longmapsto \underline{\underline{S}}(\underline{\underline{k}})$$

Now, we turn to the proof of the theorem. From Lemma 4.3, 4.5 and
 4.6,  we have
$$e_{t}=\sum_{\begin{array}{c} \underline{\underline{k}}\in \mathbb{K}^{(t+1)^{2}}\\
                  \| \underline{\underline{k}}\| = p^{t}
                \end{array} } \pi _{\underline{\underline{k}}}( \underline{a},\underline{b})=\sum_{\begin{array}{c} \underline{\underline{k}}\in \mathbb{K}^{(t+1)^{2}}\\
                 \underline{\underline{l}}(\underline{\underline{k}})\in
                 \mathbf{L}_{p}(t)
                \end{array} } \pi _{\underline{\underline{k}}}( \underline{a},\underline{b})$$
$$=\sum_{\underline{\underline{l}}\in \mathbf{L}_{p}(t)}\sum_{\begin{array}{c} \underline{\underline{k}}\in \mathbb{K}^{(t+1)^{2}}\\
\underline{\underline{l}}(\underline{\underline{k}})=(\underline{l}_{0},\ldots,\underline{l}_{m},\ldots,\underline{l}_{t})
                \end{array} }\pi _{\underline{\underline{k}}}( \underline{a},\underline{b})$$
$$=\sum_{\underline{\underline{l}}\in \mathbf{L}_{p}(t)}\sum_{(\underline{S}_{0},\ldots,\underline{S}_{m},\ldots,\underline{S}_{t})
\in \prod_{ m=0}^{ t}\mathbf{I}(m,
\underline{l}_{m})}\prod_{m=0}^{t}\prod_{j=1}^{K}\prod_{i\in
S_{m,j}}\pi_{\underline{k}(j)}(a_{i},b_{m-i})$$
$$=\sum_{\underline{\underline{l}}\in \mathbf{L}_{p}(t)}\prod_{m=0}^{t}\sum_{\underline{S}_{m}\in \mathbf{I}(m, \underline{l}_{m})}\prod_{j=1}^{K}\prod_{i\in
S_{m,j}}\pi_{\underline{k}(j)}(a_{i},b_{m-i}) \ \ \ \ \ \ \ \  \ \ \ \ \ \ \ \ \ \ \ \ \ \     $$
$$=\sum_{\underline{\underline{l}}\in \mathbf{L}_{p}(t)}\prod_{m=0}^{t}\tau_{\underline{l}_{m}}(a_{0},\ldots, a_{m}; b_{0}, \ldots, b_{m}) \ (\mod \, p). \ \ \ \ \ \ \ \  \ \ \ \ \ \ \ \ \ \ \ \ \ \    $$

$\Box$

 {\bf Corollary 4.7.} {\it  Assume that
$$a=\sum_{i=0}^{\infty}a_{i}p^{i}, b=\sum_{i=0}^{\infty}b_{i}p^{i}, ab=\sum_{i=0}^{\infty}e_{i}p^{i},
$$ with $a_{i}, b_{i},  e_{i}\in \{0,1,\ldots,p-1\}.$ Then
$e_{0} = a_{0}b_{0}\,( \mod\,p)$ and  for $t\geq 1  ,$
$$  e_{t} = \sum_{\underline{\underline{l}}=(\underline{l}_{0},\ldots,\underline{l}_{k},\ldots,\underline{l}_{t} )\in \mathbf{L}_{p}(t)}\prod_{k=0}^{t}\tau_{\underline{l}_{k}}(a_{0},\ldots,a_{k}; b_{0},\ldots, b_{k}) ( \mod  \,
                p).
$$
 In particular, if $p=2,$ we have  $e_{0} = a_{0}b_{0} (\mod\,2)$
 and  for $t\geq 1,$
$$e_{t} = \sum_{(l_{1},\ldots,l_{t})\in \mathbf{L}_{2}(t)}  \prod_{ 1 \leq k \leq t} \tau _{l_{k} }(a_{0}b_{k},a_{1}b_{k-1},\cdots ,a_{k}b_{0})  (\mod \, 2 ); $$
 if $p=3,$ we have $e_{0} = a_{0}b_{0} (\mod\,3)$  and  for
$t\geq 1,$
$$ e_{t}=\sum_{(\underline{l}_{0},\ldots,\underline{l}_{k},\ldots,\underline{l}_{t} )\in \mathbf{L}_{3}(t)}\prod_{k=0}^{t}
\sum_{\underline{S}=(S_{1}, \ldots , S_{5})\in
\mathbf{I}(k,
\underline{l}_{k})}f_{\underline{S}}(a_{0},a_{1},\ldots,a_{k};b_{0},b_{1},\ldots,b_{k})\
( \mod \, 3),$$
 where
 $$f_{\underline{S}}(a_{0},a_{1},\ldots,a_{k};b_{0},b_{1},\ldots,b_{k})=\prod_{ i_{1}\in
S_{1}}a_{i_{1}}b_{k-i_{1}}\prod_{ i_{2}\in S_{2}}
a_{i_{2}}(1-a_{i_{2}})b_{k-i_{2}}$$
$$\ \ \ \ \ \ \ \ \ \ \ \ \ \ \ \ \ \ \ \ \ \ \ \ \ \ \ \ \ \ \ \ \ \ \ \ \ \cdot\prod_{ i_{3}\in
S_{3}}a_{i_{3}}^{2}b_{k-i_{3}}(1-b_{k-i_{3}})\prod_{ i\in S_{4}\cup
S_{5}}a_{i}(1-a_{i})b_{k-i}(b_{k-i}-1).$$} $\Box$

 {\bf Remark 4.8.} (i) We can give an algorithm to determine the
set $\mathbf{L}_{2}(t).$

(ii) For $p=2,$ we once gave a rather complicated proof for the
addition formula  by simplifying the well-known recursion formulas
for the addition of Witt vectors(see [1]), but we did not know
whether the similar thing is possible for the multiplication
formula. After reading that complicated proof, Browkin  found a
simple but quite different proof for our addition formula in the
case of $p=2$ (see [2]). The present proofs, in particular those for
the results in this section, were largely inspired by the following
fact in Lucas lemma:
$$a_{t}= \left(
\begin{array}{c}
                                  A \\
                                    p^{t}\\
                                \end{array}
                              \right   )  ( \mod \,\, p ),$$
 which was first pointed in [3]. This fact was also used in [4].

{\bf Question 4.9.} How to simplify the expression of $e_{t}$
further ?

\bigskip

\section{Transformation of coefficients}

In this section, we will solve Browkin's problem. At first, we
define the required polynomials as follows.
$$f_{t}(x_{0},x_{1},\ldots,x_{t-1}):=\sum_{\lambda=0}^{ t-1 } \{\sum_{c=1}^{\frac{p-1}{2}}[(x_{\lambda}+c)^{p-1}-1]\}\prod_{\lambda < i <
t}(1-x_{i}^{p-1}),$$
$$g_{t}(y_{0},y_{1},\ldots,y_{t-1}):=\sum_{\lambda=0}^{ t-1 } \{\sum_{c=\frac{p+1}{2}}^{ p-1 }[1-  (y_{\lambda}-c  )^{p-1}]\}\prod_{\lambda < i <
t}[1-\left(y_{i}-\frac{p-1}{2}\right)^{p-1}],$$ where we also have
the convention that $\prod_{i\in \phi}=1$ for the empty set $\phi.$

{\bf Theorem 5.1.} {\it Assume that $p\geq 3$ is a prime. Let
$$A=\sum_{i}^{\infty}a_{i}p^{i}=\sum_{j}^{\infty}b_{j}p^{j} \in
\mathbb{Z}_{p}, $$ with $a_{i}\in \{0,
\pm1,\pm2,\ldots,\pm\frac{p-1}{2}\}$ and $b_{j}\in
\{0,1,\ldots,p-1\}.$ Then
$$b_{t}=a_{t}+f_{t}(a_{0},a_{1},\ldots,a_{t-1})\ (\mod  \, p). \eqno(5.1)$$
$$a_{t}=b_{t}+g_{t}(b_{0},b_{1},\ldots,b_{t-1})\ (\mod  \, p).
\eqno(5.2)$$}

{\bf Proof} \ Firstly, we prove (5.1). At first, define an index
sequence. Let $j_{0}=-1$ for the initial value. If after $k-1$
rounds ($k\geq 1$) we have $j_{k-1},$ then we go on with the
following two steps:

i) Let
$$i_{k}=\{\begin{array}{c} \infty, \ \  \mbox{if} \ \{i : j_{k-1}< i, -\frac{p-1}{2}\leq a_{i}\leq
-1\}=\phi; \ \ \ \ \ \\
 \mbox{min} \{i :  j_{k-1}< i, -\frac{p-1}{2}\leq a_{i}\leq -1\}, \
  \mbox{otherwise}.
\end{array}$$
If $i_{k}=\infty,$ then the index sequence is completed; otherwise,
go on with the next step:

ii) Let
$$j_{k}=\{\begin{array}{c} \infty, \ \  \mbox{if} \ \{i : i_{k}< i, 1 \leq a_{i}\leq \frac{p-1}{2}
\}=\phi; \ \ \ \ \ \\
 \mbox{min} \{i :  i_{k}< i, 1 \leq a_{i}\leq \frac{p-1}{2}\}, \
 \mbox{otherwise}.
\end{array}$$
If $j_{k}=\infty,$ the index sequence is completed; otherwise, go on
with the  $(k+1)$-th round.

For $k\geq 1 $ we define
$$b^{\prime}_{i}=a_{i}, j_{k-1}< i <  i_{k}, \ \   \mbox{and} \ \ b^{\prime}_{i_{k}}=p+a_{i_{k}}. \eqno(5.3)$$
$$b^{\prime}_{i}=a_{i}-1+p, i_{k}<i<j_{k},  \   \mbox{and}  \ \  b^{\prime}_{j_{k}}=a_{j_{k}}-1. \eqno(5.4)$$
It is easy to check that $0\leq b^{\prime}_{t}< p $ for any $t.$

 We
will denote
$$I_{k}=\sum_{j_{k-1}< i\leq i_{k}}a_{i}p^{i}, \ \ \ J_{k}=\sum_{i_{k}< i\leq j_{k}}a_{i}p^{i}, \ \forall k\geq 1.$$
When $i_{k}=\infty,$ from (5.3) we have
$$I_{k}=\sum_{j_{k-1}< i\leq i_{k}=\infty}a_{i}p^{i}=\sum_{j_{k-1}< i <i_{k}=\infty}a_{i}p^{i}=\sum_{j_{k-1}< i}b^{\prime}_{i}p^{i}.\eqno(5.5)$$
When $i_{k}<\infty,$ from (5.3) we have
$$I_{k}=\sum_{j_{k-1}< i\leq i_{k}}a_{i}p^{i}=\sum_{j_{k-1}< i <i_{k}}b^{\prime}_{i}p^{i}+b^{\prime}_{i_{k}}p^{i_{k}}-p^{1+i_{k}}=\sum_{j_{k-1}< i \leq i_{k}}b^{\prime}_{i}p^{i}-p^{1+i_{k}}.\eqno(5.6)$$
When $j_{k}=\infty,$ from (5.4) we have
$$-p^{1+i_{k}}+J_{k}=\sum_{i_{k}< i}(p-1)p^{i}+\sum_{i_{k}< i \leq j_{k}=\infty}a_{i}p^{i}=\sum_{i_{k}< i}(a_{i}+p-1)p^{i}=\sum_{i_{k}< i}b^{\prime}_{i}p^{i}. \eqno(5.7)$$
When $j_{k}<\infty,$ from (5.4) we have
$$-p^{1+i_{k}}+J_{k}=\sum_{i_{k}< i}(p-1)p^{i}+\sum_{i_{k}< i \leq j_{k}}a_{i}p^{i} \ \ \ \ \ \ \ \ \ \ \ \ \ \ \ \ \ \ \ \ \ \ \ \ \ \ \ \ \ \ \ \ \ \ \ \ \ \ \ \ \  $$
$$\, =\sum_{i_{k}<i<  j_{k}}(a_{i}+p-1)p^{i}+[a_{j_{k}}+\sum_{0\leq
i}(p-1)p^{i}]p^{j_{k}}$$
$$=\sum_{i_{k}< i <j_{k}}(a_{i}+p-1)p^{i}+(a_{j_{k}}-1)p^{j_{k}} \ \ \ \ \ \ \ \ \ \ \ \ \ $$
$$=\sum_{i_{k}< i \leq j_{k}}b^{\prime}_{i}p^{i}. \ \ \ \ \ \ \ \ \ \ \ \ \ \   \ \ \ \ \ \ \ \ \ \ \  \ \ \ \ \ \ \ \ \ \ \  \ \ \ \ \ \ \ \ \eqno(5.8)$$
When $j_{k}=\infty,$ from (5.6)(5.7) we have
$$I_{k}+J_{k}=\sum_{j_{k-1}<i}b^{\prime}_{i}p^{i}.\eqno(5.9)$$
When $j_{k}<\infty,$ from (5.6)(5.8) we have
$$I_{k}+J_{k}=\sum_{j_{k-1}<i\leq i_{k}}b^{\prime}_{i}p^{i}.\eqno(5.10)$$

It is easy to see that
$$A=\left \{\begin{array}{c} I_{1}+J_{1}+\cdots
+I_{k-1}+J_{k-1}+I_{k},\   \mbox{if} \ i_{k}=\infty;\\
I_{1}+J_{1}+\cdots
+I_{k}+J_{k},\   \mbox{if} \ j_{k}=\infty;\ \ \ \ \ \ \ \ \ \ \ \ \ \\
\sum_{k\geq 1}(I_{k}+J_{k}),\  \mbox{otherwise}. \ \ \ \ \ \ \ \ \ \
\ \ \ \ \ \ \ \ \ \ \ \
\end{array}\right. $$
Discussing the three cases respectively, from (5.5)-(5.10) we have
$$A=\sum_{i\geq 0}b^{\prime}_{i}p^{i}.$$

By the definition of the index sequence,  for $ k\geq1 $ clearly we
have

a) if $j_{k-1}< t \leq i_{k},$ then $0\leq a_{t-1}\leq
\frac{p-1}{2},$ and $(a_{0},a_{1},\ldots,a_{t-1})$ is not of the
form
 $(\ast, \ldots, \ast, -c,\underbrace{0,\ldots,0}_{m})$ with
$m\geq 0$ and $ 1\leq c \leq \frac{p-1}{2};$

b) if\, $i_{k}< t \leq j_{k},$ then $-\frac{p-1}{2}\leq a_{t-1}\leq
0,$ and $(a_{0},a_{1},\ldots,a_{t-1})$ is of the form $(\ast,
\ldots, \ast, -c, \underbrace{0,\ldots,0}_{m})$ with $m\geq 0$ and $
1\leq c \leq \frac{p-1}{2}.$

Hence, for $ k\geq1 $ we have  $i_{k}< t \leq j_{k} $ if and only if
$ (a_{0},a_{1},\ldots,a_{t-1})$ is  of the form $(\ast, \ldots,
\ast, -c, \underbrace{0,\ldots,0}_{m})$ with $ m\geq 0$ and $ 1\leq
c \leq \frac{p-1}{2}.$ Note that we have modulo $p$ :
$$f_{t}(a_{0},a_{1},\ldots,a_{t-1})=\{ \begin{array}{c} -1, \  \mbox{if} \ (a_{0},a_{1},\ldots,a_{t-1})=(\ast, \ldots,
\ast, -c, 0,\ldots,0),   1\leq c \leq
\frac{p-1}{2};\\
0, \  \mbox{otherwise}. \ \ \ \ \ \ \ \ \ \ \ \ \ \ \
 \ \ \ \ \ \ \ \ \ \ \ \ \ \ \ \ \ \ \ \ \ \ \ \ \ \ \ \ \ \ \ \ \ \ \ \ \ \ \ \ \ \ \ \ \ \  \end{array}
$$
So
$$a_{t}+f_{t}(a_{0},a_{1},\ldots,a_{t-1})=\{\begin{array}{c} a_{t} \ (\mod \, p), \  \mbox{if} \ j_{k-1}< t\leq i_{k},k\geq 1;\ \ \\
a_{t}-1 \ ( \mod \, p), \  \mbox{if} \ i_{k}< t\leq j_{k},k\geq
1.
\end{array}$$
Therefore, from (5.3)(5.4), we have
$$a_{t}+f_{t}(a_{0},a_{1},\ldots,a_{t-1})=b^{\prime}_{t}\ (\mod \, p). \eqno(5.11)$$
 By the uniqueness,
we have $b_{i}=b^{\prime}_{i}$ for any $i,$ so (5.1) follows from
$(5.11).$

In a similar way, we can prove $(5.2).$ Similarly, define an index
sequence. Let $j_{0}=-1$ for the initial value. If after $k$ rounds
($k\geq 1$) we have $j_{k-1},$ then we go on with the following two
steps:

i) Let
$$i_{k}=\{\begin{array}{c} \infty, \ \ \mbox{if} \ \{i : j_{k-1}< i, \frac{p-1}{2}\leq b_{i}\leq
p-1\}=\phi; \ \ \ \ \ \\
 \mbox{min}  \{i :  j_{k-1}< i, \frac{p-1}{2}\leq b_{i}\leq p-1\}, \
 \mbox{otherwise}.
\end{array}$$
If $i_{k}=\infty,$ then the index sequence is completed; otherwise,
go on with the next step:

ii) Let
$$j_{k}=\{\begin{array}{c} \infty, \ \   \mbox{if} \ \{i : i_{k}< i, 0 \leq b_{i}< \frac{p-1}{2}
\}=\phi; \ \ \ \ \ \\
 \mbox{min} \{i :  i_{k}< i, 0 \leq b_{i}< \frac{p-1}{2}\}, \
 \mbox{otherwise} .
\end{array}$$
If $j_{k}=\infty,$ the index sequence is completed; otherwise, go on
with the  $k+1$ round.

For $k\geq 1$ we define
$$a^{\prime}_{i}=b_{i}, j_{k-1}< i <  i_{k}, \ \   \mbox{and} \ \ a^{\prime}_{i_{k}}=b_{i_{k}}-p. \eqno(5.12)$$
$$a^{\prime}_{i}=b_{i}+1-p, i_{k}<i<j_{k},  \   \mbox{and} \ \  a^{\prime}_{j_{k}}=b_{j_{k}}+1. \eqno(5.13)$$
It is easy to check that $-\frac{p-1}{2}\leq a^{\prime}_{t}\leq
\frac{p-1}{2} $ for any $t.$

For $ k\geq 1,$ let
$$I_{k}=\sum_{j_{k-1}< i\leq i_{k}}b_{i}p^{i}, \ \ \ J_{k}=\sum_{i_{k}< i\leq j_{k}}b_{i}p^{i}.$$
When $i_{k}=\infty,$ from (5.12) we have
$$I_{k}=\sum_{j_{k-1}< i\leq i_{k}=\infty}b_{i}p^{i}=\sum_{j_{k-1}< i}a^{\prime}_{i}p^{i}.\eqno(5.14)$$
When $i_{k}<\infty,$ from (5.12) we have
$$I_{k}=\sum_{j_{k-1}< i\leq i_{k}}b_{i}p^{i}=\sum_{j_{k-1}< i <i_{k}}b_{i}p^{i}+b_{i_{k}}p^{i_{k}}=\sum_{j_{k-1}< i \leq i_{k}}b_{i}p^{i}+p^{1+i_{k}}.\eqno(5.15)$$
When $j_{k}=\infty,$ from (5.13) we have
$$p^{1+i_{k}}+J_{k}=-\sum_{i_{k}< i}(p-1)p^{i}+\sum_{i_{k}< i \leq j_{k}=\infty}b_{i}p^{i}=\sum_{i_{k}< i}a^{\prime}_{i}p^{i}. \eqno(5.16)$$
When $j_{k}<\infty,$ from (5.13) we have
$$p^{1+i_{k}}+J_{k}=-\sum_{i_{k}< i}(p-1)p^{i}+\sum_{i_{k}< i \leq j_{k}}b_{i}p^{i} \ \ \ \ \ \ \ \ \ \ \ \ \ \ \ \ \ \ \ \ \ \ \ \ \ \ \ \ \ \ \ \ \ \ \ \ \ \ \  $$
$$ \ \ \ \ \ \ \ \ \  \ \ \ \ \ \ =\sum_{i_{k}<i<
j_{k}}(b_{i}-p+1)p^{i}+(b_{j_{k}}+1)p^{j_{k}}-p^{1+j_{k}}-\sum_{j_{k}
< i}(p-1)p^{i}$$
$$=\sum_{i_{k}< i \leq j_{k}}a^{\prime}_{i}p^{i}. \ \ \ \ \ \ \ \ \ \ \   \ \ \ \ \ \ \ \ \ \ \  \ \ \ \ \ \ \ \ \ \ \  \ \ \ \ \ \ \ \ \ \ \ \ \ \eqno(5.17)$$
Then, similarly from (5.14)-(5.17), we have
$$A=\sum_{i\geq 0}a^{\prime}_{i}p^{i}.$$

By the definition of the index sequence, for $k\geq1$ we have:

a) if $j_{k-1}< t \leq i_{k}, $  then $0\leq b_{t-1}\leq
\frac{p-1}{2},$ and $(b_{0},b_{1},\ldots,b_{t-1})$ is not the form
of $(\ast, \ldots, \ast,
c,\underbrace{\frac{p-1}{2},\ldots,\frac{p-1}{2}}_{m})$ with $m\geq
0$ and $ \frac{p-1}{2}< c < p\, ;$

b) if \,$i_{k}< t \leq j_{k},$ then $\frac{p-1}{2}\leq b_{t-1} < p\,
,$ and $(b_{0},b_{1},\ldots,b_{t-1})$ is the form of $(\ast, \ldots,
\ast, c,\underbrace{\frac{p-1}{2},\ldots,\frac{p-1}{2}}_{m})$ with
$m\geq 0$ and $ \frac{p-1}{2}< c <  p\, .$

Therefore, for $k\geq1$ we have that $i_{k}< t \leq j_{k} $ if and
only if $ (b_{0},b_{1},\ldots,b_{t-1})$ is the form of $(\ast,
\ldots, \ast,
c,\underbrace{\frac{p-1}{2},\ldots,\frac{p-1}{2}}_{m})$ with $m\geq
0$ and $ \frac{p-1}{2}< c <  p.$ Note that we have modulo $p$ :
$$g_{t}(b_{0},b_{1},\ldots,b_{t-1})=\{ \begin{array}{c} 1, \  \mbox{if} \ (b_{0},b_{1},\ldots,b_{t-1})=(\ast, \ldots, \ast,
c,\frac{p-1}{2},\ldots,\frac{p-1}{2}),  \frac{p-1}{2}< c <  p;\\
0 , \  \mbox{otherwise} . \ \ \ \ \ \ \ \ \ \ \ \ \ \ \
 \ \ \ \ \ \ \ \ \ \ \ \ \ \ \ \ \ \ \ \ \ \ \ \ \ \ \ \ \ \ \ \ \ \ \ \ \ \ \ \ \ \ \ \ \ \ \ \ \  \end{array}
$$
So
$$b_{t}+g_{t}(b_{0},b_{1},\ldots,b_{t-1})=\{\begin{array}{c} b_{t}+1 \ (\mod  p), \  \mbox{if} \ j_{k-1}< t\leq i_{k},k\geq 1;\ \ \\
b_{t} \ (\mod p), \  \mbox{if} \ i_{k}< t\leq j_{k},k\geq 1. \
\ \ \ \ \ \ \ \ \ \
\end{array}$$
Hence
$$b_{t}+g_{t}(b_{0},b_{1},\ldots,b_{t-1})=a^{\prime}_{t}\ (\mod p). \eqno (5.18)$$
As above, by uniqueness we know that (5.2) follows from (5.18).
$\Box$

{\bf An alternative proof} After read the previous version of this
paper, Browkin gave an alternative proof for Theorem 5.1. Now, we
only give a  sketch of his proof of the equality (5.1).

Let $\sum_{i=0}^{\infty}a_{i}p^{i}=\sum_{i=0}^{\infty}b_{i}p^{i},$
where $a_{i}\in \{0, \pm1,\pm2,\ldots,\pm\frac{p-1}{2}\}, b_{i}\in
\{0, 1,  \ldots \\ , p-1\}.$ For $k\geq 0$ denote
$$A_{k}:=\sum_{i=0}^{k}a_{i}p^{i},\ \ \ \  B_{k}:=\sum_{i=0}^{k}b_{i}p^{i}.$$
Clearly, for any $k\geq 0$, $A_{k}, B_{k}$ satisfy $A_{k}\equiv
B_{k} (\mod p^{k+1}).$ We have
$$ \mid A_{k}\mid< p^{k+1}
\ \  \mbox{and} \ \ 0\leq B_{k}< p^{k+1}.\eqno (\ast)$$ In fact, we
have
$$\mid A_{k}\mid\leq \sum_{i=0}^{k}\mid a_{i}\mid p^{i}\leq \frac{p-1}{2}\sum_{i=0}^{k} p^{i}=\frac{1}{2}(p^{k+1}-1)< p^{k+1}$$
and
$$0\leq  B_{k}= \sum_{i=0}^{k} b_{i} p^{i}\leq (p-1)\sum_{i=0}^{k} p^{i}=p^{k+1}-1< p^{k+1}.$$
From $(\ast),$ it follows that
$$-p^{k+1}< -A_{k}\leq B_{k}-A_{k}\leq B_{k}+\mid A_{k}\mid< p^{k+1},$$
so we have $B_{k}-A_{k}=0$ or $p^{k+1}.$ More precisely
$$B_{k}=A_{k} \  \mbox{if} \ A_{k}\geq 0; \ \  B_{k}=A_{k}+p^{k+1}\  \mbox{if} \ A_{k}< 0. \eqno (\ast\ast)$$
From this, we know that $b_{0}\equiv a_{0}(\mod p).$ Now, we
determine $b_{k} (\mod p)$ for $k\geq 1.$

i) Assume that $A_{k-1}\geq 0.$ Then from $(\ast)$ we have
$A_{k-1}=B_{k-1}.$ If $A_{k}\geq 0,$ then $A_{k}=B_{k}$ similarly,
so
$$A_{k-1}+a_{k}p^{k}=A_{k}=B_{k}=B_{k-1}+b_{k}p^{k},$$
 therefore
$b_{k}=a_{k};$ if $A_{k}<0,$ then  by $(\ast\ast)$ we have
$B_{k}=A_{k}+p^{k+1},$ and so
$$B_{k-1}+b_{k}p^{k}=B_{k}=A_{k}+p^{k+1}=A_{k-1}+a_{k}p^{k}+p^{k+1},$$
which implies $b_{k}=a_{k}+p\, .$

ii) Assume that $A_{k-1}< 0.$ If $A_{k}\geq 0,$ then from
$(\ast\ast)$ we get
$$A_{k-1}+p^{k}+b_{k}p^{k}=B_{k-1}+b_{k}p^{k}=B_{k}=A_{k}=A_{k-1}+a_{k}p^{k},$$
therefore $b_{k}=a_{k}-1; $ if $A_{k}< 0,$ then from $(\ast\ast)$ we
get
$$A_{k-1}+p^{k}+b_{k}p^{k}=B_{k-1}+b_{k}p^{k}=B_{k}=A_{k}+p^{k+1}=A_{k-1}+a_{k}p^{k}+p^{k+1},$$
therefore $b_{k}=a_{k}+p-1\equiv a_{k}-1 (\mod p)$.

Thus we have proved:
$$b_{k}-a_{k}\equiv
\{\begin{array}{c} -1 (\mod p), \  \mbox{if}\ A_{k-1} < 0\, ; \\
0 \ \ (\mod p), \   \mbox{otherwise} .\end{array}$$ Now we
express these conditions by means of polynomials.

Let
$$A_{k-1}=\sum_{i=0}^{k-1}a_{i}p^{i}, \  \mbox{where}  \ a_{k}=a_{k-1}=\ldots =a_{m+1}=0, a_{m}\neq 0,$$
for some $m, 0\leq m\leq k.$ From $A_{k-1}=A_{m}=A_{m-1}+a_{m}p^{m}$
and $\mid A_{m-1}\mid < p^{m}$ we conclude that $A_{k-1}< 0$ if and
only if $a_{m}< 0,$ which is equivalent to $a_{m}\in \{-1,-2,
\ldots, -\frac{p-1}{2}\}.$ So we get
$$b_{k}-a_{k}\equiv
\{\begin{array}{c} -1 (\mod p), \  \mbox{if} \
(a_{0},a_{1},\ldots,a_{k-1})=(\ast, \ldots, \ast, -c, 0,\ldots,0); \
\\ 0 \ \ (\mod p), \  \mbox{otherwise} , \ \ \ \ \ \ \ \ \ \ \
\ \ \ \ \ \ \ \ \ \ \ \ \ \ \ \ \ \ \ \ \ \ \ \ \ \ \ \ \ \ \ \ \ \
\end{array}$$
where $1\leq c \leq \frac{p-1}{2}.$ From the proof of Theorem 5.1,
we know that $f_{k}(a_{0},a_{1},\ldots,a_{k-1})$ has the same
property as $b_{k}-a_{k},$ so we have
$$b_{k}=a_{k}+f_{k}(a_{0},a_{1},\ldots,a_{k-1})\ (\mod
\, p). $$ $\Box$

{\bf Corollary 5.2.} {\it  Let
$$A=\sum_{i}^{\infty}a_{i}3^{i}=\sum_{j}^{\infty}b_{j}3^{j} \in
\mathbb{Z}_{3}, $$ with $a_{i}\in \{0, \pm1\}$ and $b_{j}\in \{0,1,2
\}.$ Then
$$b_{t}=a_{t}+\sum_{0\leq \lambda < t } a_{\lambda}(a_{\lambda}-1)\prod_{\lambda < i <
t}(1-a_{i}^{2})\ (\mod \, 3).$$
$$a_{t}=b_{t}+\sum_{0\leq \lambda < t } b_{\lambda}(1-b_{\lambda})\prod_{\lambda < i <
t}b_{i}(2-b_{i})\ (\mod \, 3).$$} $\Box$

 We can also give the formulas of the sum and the
multiplication of $p$-adic integers with respect to the numerically
least residue system $\{0,\pm1,\pm 2,\ldots, \pm\frac{p-1}{2}\}$.
Define
$$a_{t}^{\vee}:=a_{t}+ \sum_{\lambda=0}^{ t-1 } \{\sum_{c=1}^{\frac{p-1}{2}}[(a_{\lambda}+c)^{p-1}-1]\}\prod_{\lambda < i <
t}(1-a_{i}^{p-1}),$$
$$b_{t}^{\wedge}:=b_{t}+\sum_{\lambda=0}^{ t-1 } \{\sum_{c=\frac{p+1}{2}}^{ p-1 }[1-  (b_{\lambda}-c  )^{p-1}]\}\prod_{\lambda < i <
t}[1-\left(b_{i}-\frac{p-1}{2}\right)^{p-1}],$$ where $a_{i}\in \{0,
\pm1,\pm2,\ldots,\pm\frac{p-1}{2}\}$ and $b_{j}\in
\{0,1,\ldots,p-1\}.$

{\bf Theorem 5.3. } {\it Let $p$ be an odd prime.  Assume that
$$a=\sum_{i=0}^{\infty}a_{i}p^{i}, b=\sum_{i=0}^{\infty}b_{i}p^{i}, -a=\sum_{i=0}^{\infty}d_{i}p^{i}, a+b=\sum_{i=0}^{\infty}c_{i}p^{i}\in
\mathbb{Z}_{p},ab=\sum_{i=0}^{\infty}e_{i}p^{i},
$$
 with $a_{i}, b_{i}, c_{i}, d_{i}\in \{0,\pm1,\pm 2,\ldots,
\pm\frac{p-1}{2}\}.$ Then}

(i) {\it $c_{0} = a_{0}+b_{0} (\mod\,p)$
 and  for $t\geq 1,$
$$c_{t}=a_{t}+b^{\vee}_{t}+\sum_{i=0}^{t-1}\left(\sum_{j=1}^{p-1}\left(\begin{array}{c} \frac{p-1}{2}+a_{i} \\ j\\
 \end{array}
  \right)\left(\begin{array}{c} b_{i}^{\vee} \\ p-j\\
 \end{array}
  \right)\right )\prod_{j=i+1} ^{t-1}\left(\begin{array}{c} \frac{p-1}{2}+a_{j}+b^{\vee}_{j} \\ p-1\\
 \end{array}
  \right)(\mod \ p).$$
  In particular,  if $p=3,$ then $c_{0} = a_{0}+b_{0}^{\vee}\, (\mod\,3)$  and for $t\geq
1,$
$$c_{t}=a_{t}+b^{\vee}_{t}-\sum_{i=0}^{t-1}[(a_{i}+1)(a_{i}+b_{i}^{\vee}-1)b_{i}^{\vee}]\prod_{j=i+1} ^{t-1}\left(\begin{array}{c} a_{j}+b_{j}^{\vee}+1 \\ 2\\
 \end{array}
  \right)\, (\mod\,3).$$}

(ii) {\it $d_{0}=-a_{0}^{\vee} (\mod \,p)$  and for $t\geq1$
$$d_{t}=-a_{t}^{\vee}-1+ \prod_{i=0}^{t-1}(1-{a_{i}^{\vee}}^{p-1})( \mod  \,
p).$$ In particular, if $p=3,$ then $d_{0}=-a_{0}^{\vee}
(\mod \,3)$ and for $t\geq1$
$$d_{t}=-a_{t}^{\vee}-1+ \prod_{i=0}^{t-1}(1-{a_{i}^{\vee}}^{2})(\mod  \,
3).$$}

(iii) {\it $e_{0} = (a_{0}^{\vee}b_{0}^{\vee})^{\wedge} ( \mod \,p)$
and for $t\geq 1 ,$
$$  e_{t} =\left( \sum_{\underline{\underline{l}}=(\underline{l}_{0},\ldots,\underline{l}_{k},\ldots,\underline{l}_{p} )\in \mathbf{L}_{p}(t)}\prod_{k=0}^{t}\tau_{\underline{l}_{k}}(a_{0}^{\vee},\ldots,a_{k}^{\vee}; b_{0}^{\vee},\ldots, b_{k}^{\vee})\right)^{\wedge} ( \mod  \,
                p).
$$}

{\bf Proof} (i) From Theorem 5.1, we have
$$a+b=\sum_{i=0}^{\infty}a_{i}p^{i}+\sum_{i=0}^{\infty}b^{\vee}_{i}p^{i}=\sum_{i=0}^{\infty}\left(\frac{p-1}{2}+a_{t-1}\right)p^{i}+\sum_{i=0}^{\infty}b^{\vee}_{i}p^{i}-\sum_{i=0}^{\infty}\left(\frac{p-1}{2}\right)p^{i}.$$
Note that $\frac{p-1}{2}+a_{t-1}, b ^{\vee}_{i} \in
\{0,1,\ldots,p-1\}. $ Let
$$\sum_{i=0}^{\infty}\left(\frac{p-1}{2}+a_{t-1}\right)p^{i}+\sum_{i=0}^{\infty}b^{\vee}_{i}p^{i}=\sum_{i=0}^{\infty}c_{i}^{\prime}p^{i},
\ c_{i}^{\prime}\in \{0,1,\ldots,p-1\}.$$ Then by Theorem 6.1 we
have
$$c_{t}^{\prime}=\frac{p-1}{2}+a_{t}+b^{\vee}_{t}+\sum_{i=1}^{p-1}\left(\begin{array}{c} \frac{p-1}{2}+a_{t-1} \\ i \\
 \end{array}
  \right)\left(\begin{array}{c} b^{\vee}_{t-1} \\ p-i \\
 \end{array}
  \right)$$
  $$+\sum_{i=0}^{t-2}\left(\sum_{j=1}^{p-1}\left(\begin{array}{c} \frac{p-1}{2}+a_{i} \\ j\\
 \end{array}
  \right)\left(\begin{array}{c} b_{i}^{\vee} \\ p-j\\
 \end{array}
  \right)\right )\prod_{j=i+1} ^{t-1}\left(\begin{array}{c} \frac{p-1}{2}+a_{j}+b^{\vee}_{j} \\ p-1\\
 \end{array}
  \right)(\mod \ p).$$
Clearly $c_{t}=c_{t}^{\prime}-\frac{p-1}{2}.$

(ii) It follows from Theorem  5.1 and Theorem 3.1.

(iii) It follows from Theorem 5.1, Corollary 2.4 and Corollary 4.7.
$\Box$

\bigskip

\section{Applications to Witt vectors}

 Now, we apply the above results to
$(\mathbf{W}(\mathbb{F}_{p}),\dot{+}, \dot{\times})$, the ring of
Witt vectors with coefficients in $\mathbb{F}_{p}.$ Let $ \dot{-}$
denote the minus of Witt vectors.

{\bf Theorem 6.1. } {\it Let $a=(a_{0}, a_{1},\ldots,a_{n},\ldots),
b=(b_{0}, b_{1},\ldots,b_{n},\ldots)\in \mathbf{W}(\mathbb{F}_{2})$.
If in $ \mathbf{W}(\mathbb{F}_{2})$
 $$a\dot{+}b=(c_{0}, c_{1},\ldots,c_{n},\ldots),$$
$$ \dot{-}a=(d_{0}, d_{1},\ldots,d_{n},\ldots),$$
$$a\dot{\times}b=(e_{0}, e_{1},\ldots,e_{n},\ldots),$$
 then  in $\mathbb{F}_{2}$ we have}

(i) {\it $c_{0}=a_{0}+b_{0}$ and for} $t\geq1,$
 $$c_{t} = a_{t}+b_{t}+\sum_{i=0}^{t-1}a_{i}b_{i}\prod_{j=i+1}^{t-1}(a_{j}+b_{j}).$$

(ii) {\it  $d_{0}=a_{0},$ and for} $t\geq1,$
$$ d_{t} = a_{t} +  1+ \prod_{i=0}^{t-1}(1+a_{i}).$$

(iii) {\it $e_{0}=a_{0}b_{0},$ and for} $t\geq 1,$
$$e_{t} = \sum_{(l_{1},\ldots,l_{t})\in \mathbf{L}_{2}(t) }  \prod_{ 1 \leq k \leq t} \tau _{l_{k} }(a_{0}b_{k},a_{1}b_{k-1},\cdots ,a_{k}b_{0}). $$

{\bf Proof} It follows from Corollary 2.4 and 4.7. $\Box$

When $p=3$, $a_{t}^{\vee}$ and $ b_{t}^{\wedge}$ become
$$a_{t}^{\vee}= a_{t}+\sum_{0\leq \lambda < t } a_{\lambda}(a_{\lambda}-1)\prod_{\lambda < i <
t}(1-a_{i}^{2}),$$
$$b_{t}^{\wedge}=b_{t}+\sum_{0\leq \lambda < t } b_{\lambda}(1-b_{\lambda})\prod_{\lambda < i <
t}b_{i}(2-b_{i}) $$ with $a_{i}\in \{0, \pm1\}$ and $b_{j}\in
\{0,1,2\},$ and then we have:

{\bf Theorem 6.2. } {\it Let $a=(a_{0}, a_{1},\ldots,a_{n},\ldots),
b=(b_{0}, b_{1},\ldots,b_{n},\ldots)\in \mathbf{W}(\mathbb{F}_{3})$,
 If in $
\mathbf{W}(\mathbb{F}_{3})$
 $$a\dot{+}b=(c_{0}, c_{1},\ldots,c_{n},\ldots),$$
 $$ \dot{-}a=(d_{0}, d_{1},\ldots,d_{n},\ldots),$$
$$a\dot{\times}b=(e_{0}, e_{1},\ldots,e_{n},\ldots),$$
 then in $\mathbb{F}_{3}$ we have}

 (i) {\it $c_{0} = a_{0}+b^{\vee}_{0}$ and  for} $t\geq 1,$
$$c_{t}=a_{t}+b^{\vee}_{t}-\sum_{i=0}^{t-1}[(a_{i}+1)(a_{i}+b_{i}^{\vee}-1)b_{i}^{\vee}]\prod_{j=i+1} ^{t-1}\left(\begin{array}{c} a_{j}+b_{j}^{\vee}+1 \\ 2\\
 \end{array}
  \right).$$

(ii) {\it  $d_{0}=-a_{0}^{\vee}$  and for $t\geq1$
$$d_{t}=-a_{t}^{\vee}-1+ \prod_{i=0}^{t-1}(1-{a_{i}^{\vee}}^{2}).$$}

  (iii) $e_{0} = (a_{0}^{\vee}b_{0}^{\vee})^{\wedge} $ {\it and  for} $t\geq 1  ,$
$$ e_{t}=\left (\sum_{(\underline{l}_{0},\ldots,\underline{l}_{k},\ldots,\underline{l}_{t} )\in \mathbf{L}_{3}(t)}\prod_{k=0}^{t}
\sum_{\underline{S}=(S_{1}, \ldots , S_{5})\in
\mathbf{I}(k,
\underline{l}_{k})}f^{\vee}_{\underline{S}}(a_{0},a_{1},\ldots,a_{k};b_{0},b_{1},\ldots,b_{k})\right)^{\wedge},$$
 where
 $$f^{\vee}_{\underline{S}}( a_{0},a_{1},\ldots,a_{k};b_{0},b_{1},\ldots,b_{k})=\prod_{
i_{1}\in S_{1}}a_{i_{1}}^{\vee}b_{k-i_{1}}^{\vee}\prod_{ i_{2}\in
S_{2}} a_{i_{2}}^{\vee}(1-a_{i_{2}}^{\vee})b_{k-i_{2}}^{\vee}\ \ \ \
\ \ \ \ \ \ \ \
$$
$$\ \ \ \ \ \ \ \ \ \ \ \ \ \ \ \ \ \ \ \ \ \ \ \ \ \ \ \ \ \  \cdot\prod_{
i_{3}\in
S_{3}}{a_{i_{3}}^{\vee}}^{2}b^{\vee}_{k-i_{3}}(1-b^{\vee}_{k-i_{3}})
 \cdot\prod_{ i\in
S_{4}\cup
S_{5}}{a^{\vee}_{i}}^{2}(1-a^{\vee}_{i})b^{\vee}_{k-i}(b^{\vee}_{k-i}-1)$$

{\bf Proof} \ It follows from Corollary 2.4, Corollary 4.7 and
Theorem 5.3 (See [1]). $\Box$

{\bf Remark 6.3.}  (i) We can also write out for Witt vectors the
results corresponding Corollary 2.5 and  2.6.

(ii) The formulas given in Theorem 6.2 in particular for  $e_{t}$
are really terribly complicated, but they are patterns.

{\bf Question 6.4.} Can we give similar formulas for $
\mathbf{W}(\mathbb{F}_{p})$ for a prime $p> 3$ ?

\bigskip

{\bf Acknowledgment } We are grateful
 to J. Browkin for many helpful suggestions, in particular, for his  showing us the problems.

\bigskip

\bigskip

\bigskip

{\begin{center} {\bf References}
\end{center}}

\bigskip

[1] J. P. Serre, {\it Local Fields}, Springer-Verlag, New York
Heidelberg Berlin, 1979.

[2] J. Browkin, {\it The sum of dyadic numbers}, preprint.

[3] F.J.Macwilliam and N.J.A.Sloane, {\it The Theory of
Error-Correcting Codes}, North-Holland Publishing Company, 1977

[4] Bao Li and Zongduo Dai, {\it A general result and a new lower
bound of linear complexity for binary sequences derived from
sequences over $\mathbb{Z}_{2^e}$}, preprint.

\bigskip

\bigskip

\bigskip

\end{document}